\documentclass [11pt, a4paper, twoside]{article}

\usepackage {a4, indentfirst, latexsym, theorem, QED, amssymb, epsfig}

\renewcommand{\qedsymbol}{\quad\Box}
\renewcommand{\TheWordProof}{\bf Proof.}

\newcommand\nothing[1]{\relax}

\theoremstyle{change}
\theorembodyfont{\it}
\newtheorem{theorem}{Theorem.}[section]
\newtheorem{proposition}[theorem]{Proposition.}
\newtheorem{lemma}[theorem]{Lemma.}

{\theorembodyfont{\upshape}
\newtheorem{definition}[theorem]{Definition.}
}

\newcommand{\Vis}{\mathrm{Vis }}
\newcommand{\Invis}{\mathrm{Inv }}
\newcommand{\Back}{\mathrm{Back }}
\newcommand{\Front}{\mathrm{Front }}
\newcommand{\Low}{\mathrm{Low }}
\newcommand{\Up}{\mathrm{Up }}
\newcommand{\bR}{{\mathbb R}}
\newcommand{\bZ}{{\mathbb Z}}
\newcommand{\bC}{{\mathbb C}}
\newcommand{\iso}{\cong}
\newcommand{\ie}{{\it i.e.}}

\date{}

\begin{document}

\title{On a theorem of Brion}
\author{Thomas H\"uttemann}
\maketitle

\centerline {\it University of Leicester, Department of Mathematics}
\centerline {\it University Road, Leicester LE1~7RH, England (UK)}
\centerline {e-mail: \texttt{th68@mcs.le.ac.uk}}

\vglue 2\bigskipamount \hrule \medskip

{\footnotesize \smallskip We give an elementary geometric re-proof of
  a formula discovered by \textsc{Brion\/} as well as two variants
  thereof. A subset~$K$ of $\bR^n$ gives rise to a formal
  \textsc{Laurent\/} series with monomials corresponding to lattice
  points in~$K$. Under suitable hypotheses, these series represent
  rational functions~$\sigma(K)$. We will prove formul\ae{} relating
  the rational function $\sigma(P)$ of a lattice polytope~$P$ to the
  sum of rational functions corresponding to the supporting cones
  subtended at the vertices of~$P$. The exposition should be suitable
  for everyone with a little background in topology.


\noindent
{\it Keywords:\/} Polytope, cone, lattice point, generating function,
lattice point enumerator, Brion's formula\hfill\today}
\medskip \hrule \vglue 2\bigskipamount

\section{Brion's formula}

The goal of this note is to exhibit a geometric proof of an astonishing
formula discovered by \textsc{Brion}, relating the lattice point
enumerator of a rational polytope to the lattice points enumerators of
supporting cones subtended at its vertices. Roughly speaking, the
theorem is about the surprising fact that in a certain sum of rational
functions which are all given by infinite \textsc{Laurent\/} series,
there is enough cancellation so that only finitely many terms survive:
The sum collapses to a \textsc{Laurent\/} polynomial.

The argument is based on systematic usage of \textsc{Euler}
characteristics of visibility complexes. The method of proof will
readily yield two variants of \textsc{Brion}'s formula as well. These
notes are intended as an easily accessible introduction for
non-experts with some topological background.

We start with a $1$-dimensional example to demonstrate the
cancellation. The first series we consider is the well-known geometric
series $\sum_{j=0}^\infty x^j$. For $|x| < 1$ this series converges
to the rational function $f_1 (x) = 1/(1-x)$. The second series is a
variant of the geometric series in $x^{-1}$, namely
$\sum_{j=-\infty}^2 x^j = x^2 \sum_{j=-\infty}^0 x^j$. For
$|x^{-1}|<1$ this series converges to $f_2(x) = x^2 / (1-x^{-1})$. The two
series have no common domain of convergence; we can, however, add the
rational functions they represent and obtain
\[f_1(x) + f_2(x) = {1 \over 1-x} + {x^2 \over 1-x^{-1}} = {1 \over
  1-x} + {x^3 \over x-1} = {1 - x^3 \over 1-x} = 1+x+x^2 \ ,\]
a polynomial with three terms only. Note that this happens only on the
level of rational functions; adding the power series yields a
(non-convergent, formal) power series with infinitely many terms.--The
interested reader might want to check the paper
\cite{Beck:Brion_formulae} which contains a careful exposition of a
$2$-dimensional example.

\medskip

To formulate the main theorem, and to link the example to geometry, we
have to introduce some notation first. Given a subset $K \subseteq
\bR^n$ and a vector $b \in \bR^n$ we define $b+K$ as the set of points
of~$K$ shifted by the vector~$b$:
\[b+K = \{ b + x \,|\, x \in K \} \ .\]
We associate to each subset~$K$ of~$\bR^n$ a formal \textsc{Laurent}
series~$\mathcal{S}(K)$ with complex coefficients in $n$
indeterminates as follows. We write $\bC[[x_1^{\pm 1}, x_2^{\pm 1},
\ldots, x_n^{\pm 1}]]$ for the set of \textsc{Laurent\/} series; it is
a module over the ring $\bC [x_1^{\pm 1}, x_2^{\pm 1}, \ldots,
x_n^{\pm 1}]$ of \textsc{Laurent\/} polynomials.

For a given vector $\mathbf{a} = (a_1, a_2, \ldots, a_n) \in \bZ^n$ we
write $x^\mathbf{a}$ for the product $x_1^{a_1} x_2^{a_2} \ldots
x_n^{a_n}$.

\begin{definition}
  For a subset $K \subseteq \bR^n$ we define the formal \textsc{Laurent}
  series
  \[\mathcal{S}(K) = \sum_{\mathbf{a} \in \bZ^n \cap K} x^{\mathbf{a}}
  \ \in \bC[[x_1^{\pm 1}, x_2^{\pm 1}, \ldots, x_n^{\pm 1}]] \ .\]
\end{definition}

A straightforward calculation shows $\mathcal{S}(\mathbf{b}+K) =
x^\mathbf{b} \mathcal{S}(K)$ for any $\mathbf{b} \in \bZ^n$.

In favourable cases the series $\mathcal{S}(K)$ represents a rational
function which we will denote $\sigma (K) \in \bC (x_1, x_2, \ldots,
x_n)$. For example, if $K = (-\infty, 2] \subset \bR$, then
$\mathcal{S}(K) = \sum_{j=-\infty}^2 x^j$, so $\sigma(K) = x^2/(1-x^{-1})$.

\medskip

Let $P$ denote a polytope (the convex hull of a finite set of points)
in~$\bR^n$. We assume throughout that $P$ has non-empty interior, \ie,
that $P$ is of dimension~$n$.
Given a vertex $v$ of~$P$ we define the {\it barrier cone\/} $C_v$
of~$P$ at~$v$ as the set of finite linear combination with
non-negative real coefficients spanned by the set $-v + P$.
This is a cone based at the origin of the coordinate system, having
the origin as a vertex. It is the smallest such cone containing the
translate $-v+P$ of~$P$.

Since $C_v$ is a pointed cone, the associated \textsc{Laurent\/}
series represents a rational function. We can thus formulate the
following result (where $-P = \{-x \,|\, x \in P\}$ in
Equation~(\ref{eq:brion_intP})):

\begin{theorem}
  \label{thm:brion}
  Suppose $P$ is a polytope such that all its faces admit rational
  normal vectors (this happens, for example, if~$P$ has vertices
  in~$\bZ^n$).  Then there are equalities of rational functions
  \begin{eqnarray}
    \label{eq:brion_P}
    \sum_{v {\rm\ vertex\ of\ } P} \sigma (v+C_v) &=& \sigma (P) \ ,\\
    \noalign{\vglue 1 em}
    \label{eq:brion_1}
    \sum_{v {\rm\ vertex\ of\ } P} \sigma (C_v) &=& 1 \ ,\\
    \noalign{\vglue 1 em}
    \label{eq:brion_intP}
    \sum_{v {\rm\ vertex\ of\ } P} \sigma (-v+C_v) &=& (-1)^n \sigma
    (\mathrm{int}\,(-P)) \ .
  \end{eqnarray}
\end{theorem}

A version of this theorem appears as Proposition~3.1
in~\cite{Brion-Vergne}. We give a more geometric proof, working out
how visibility subcomplexes of polytopes enter the picture.

\medskip

The example above shows how Equation~(\ref{eq:brion_P}) works for
$P = [0,2] \subset \bR$. The vertices of~$P$ are $v=0$ and
$v=2$, the respective barrier cones are
\[C_0 = [0, \infty) \qquad \hbox{and} \qquad C_2 = (-\infty, 0] \ ,\]
so $0 + C_0 = C_0$ and $2 + C_2 = (-\infty, 2]$. For
Equation~(\ref{eq:brion_1}), since $\mathcal{S}(C_0) =
\sum_{j=0}^\infty x^j$ we have $\sigma (C_0) = 1/(1-x)$, and similarly
$\sigma (C_2) = 1/(1-x^{-1})$, so the Theorem predicts correctly that
\[\sigma(C_0) + \sigma(C_2) = {1 \over 1-x} + {1 \over 1-x^{-1}} =
{1-x \over 1-x} = 1 \ .\]
Finally, we consider Equation~(\ref{eq:brion_intP}). We have $-P =
[-2,0]$, so the only integral point in the interior of $-P$ is $-1$,
and the right-hand side of~(\ref{eq:brion_intP}) is the single
term~$-x^{-1}$. On the left, we have $\sigma (-0+C_0) = 1/(1-x)$ as
before, and $\sigma(-2 + C_2) = x^{-2}/ (1-x^{-2})$, and indeed
\[\sigma(-0+C_0) + \sigma (-2+C_2) = {1 \over 1-x} + {x^{-2} \over
  1-x^{-1}} = {1-x^{-1} \over 1-x} = x^{-1} \cdot {x-1 \over 1-x} =
-x^{-1} \ .\]

\medskip

We will prove Theorem~\ref{thm:brion} in
\S\S\ref{sec:barr-cones}--\ref{sec:brions-formula}.
Equation~(\ref{eq:brion_P}) of the Theorem is the original version of
\textsc{Brion}'s formula \cite[Theorem~2.2]{Brion-latticepoints}
\cite[Theorem~2.1~(ii)]{Brion-survey}. 

The paper is inspired by \textsc{Beck}, \textsc{Haase} and
\textsc{Sottile} \cite{Beck:Brion_formulae} who gave a new, elementary
proof of \textsc{Brion}'s formula.  The approach taken in this note is
a rather straightforward elaboration: It replaces the elegant but
delicate combinatorics of \cite{Beck:Brion_formulae} with a geometric
analysis of visibility subcomplexes
(\S\ref{subsec:VisibilityComplexes}) of a polytope. From a
topologist's point of view this makes the proof more transparent,
while still avoiding the elaborate machinery of toric algebraic
geometry used in the original proof \cite{Brion-latticepoints}.

The basic strategy of proof is to establish an identity of formal
\textsc{Laurent }power series first (Theorem~\ref{thm:brion_powerseries}),
then pass to rational functions (Theorem~\ref{thm:brion}). This is explained
in detail in \cite{Beck:Brion_formulae}, but we will recall the relevant
arguments for the convenience of the reader.

\section{Polytopal complexes}

\label{sec:polytopal-complexes}

A {\it polytope\/} $P$ is the convex hull of a non-empty finite set of
points in~$\bR^n$. A face of~$P$ is the intersection of~$P$ with some
supporting hyperplane; as a matter of convention, we also have the
improper faces $F = P$ and $F=\emptyset$. See \cite{Ewald-CCAG} and
\cite{Ziegler-Polytopes} for more on polytopes and their faces.

\begin{definition}
A non-empty finite collection~$K$ of non-empty polytopes
in some~$\bR^n$ is called a
{\it polytopal complex\/} if the following conditions
are satisfied:
\begin{enumerate}
        \item If $F \in K$ and $G$ is a non-empty face of~$F$, then $G \in K$.
        \item For all $F, G \in K$, the intersection $F \cap G$
        is a (possibly empty) face of both~$F$ and~$G$.
\end{enumerate}
A subset $L \subseteq K$ of a polytopal complex is called
an {\it order filter\/} if for all $F \in L$ and $G \in K$
with $F$ a face of~$G$, we have $G \in L$.
A subset $L \subseteq K$ of a polytopal complex is called
a {\it subcomplex\/} of~$K$ if $L$~is a polytopal complex.
\end{definition}

Important examples of polytopal complexes are the complex~$F(P)_0$ of non-empty
faces of a polytope~$P$, and its subcomplex~$F(P)^1_0$ of non-empty
proper faces of~$P$ (sometimes called boundary complex of~$P$).

The intersection of two subcomplexes, if non-empty, is a subcomplex. The
(set-theoretic) complement of a subcomplex is an order filter.

\begin{definition}
Suppose $K$~is a polytopal complex, and $L$~is a non-empty subset of~$K$. We
call $ |L| := \bigcup_{F \in L} F $ the {\it realisation\/} or the
{\it underlying space\/}
of~$L$.
\end{definition}

If~$P$~is an $n$-dimensional polytope, we have homeomorphisms
$|F(P)_0| = P \iso B^n$ and $|F(P)_0^1| = \partial P \iso S^{n-1}$.

\begin{definition}
  Let $L$ be a non-empty subset of the polytopal complex~$K$. The {\it
  \textsc{Euler\/} characteristic\/} $\chi(L)$ of~$L$ is defined by
  \[\chi(L) = \sum_{A \in L} (-1)^{\dim (A)} \ .\]
\end{definition}

If $L$ is a subcomplex of~$K$, then $\chi(L)$ agrees with the \textsc{Euler\/}
characteristic of~$|L|$ as defined in algebraic topology. In particular,
$\chi(F(P)_0) = \chi(P) = 1$ and $\chi (F(P)_0^1) = \chi (\partial P) = 1 +
(-1)^{\dim (P)}$ for any polytope~$P$.

\begin{lemma}
  The \textsc{Euler\/} characteristic is additive: For a polytopal
  complex~$K$ and a non-empty proper subset $L \subset K$ we have
  \[\chi(K) = \chi (L) + \chi (K \setminus L) \ . \qed\]
\end{lemma}


\section{Visibility subcomplexes of a polytope}
\label{subsec:VisibilityComplexes}

Understanding visibility subcomplexes of polytopes is the key to our
approach to \textsc{Brion}'s theorem. The notions of visible, back and
lower faces are defined, and we indicate a proof that these
subcomplexes are balls in the boundary sphere of~$P$. In particular,
these complexes are contractible and have \textsc{Euler\/}
characteristic~$1$.---We assume throughout that $P \subset \bR^n$ is a
polytope with $\mathrm{int}(P) \not= \emptyset$.

\subsection*{Visible and invisible faces}
\label{SecVisInvis}

\begin{definition}
  A face $ F \in F(P)_0^1 $ is called {\it visible\/} from the point $x \in
  \bR^n \setminus P$ if $ [p,x] \cap P = \{p\} $ for all $ p \in F $. (Here
  $[p,x]$ denotes the line segment between~$p$ and~$x$.) Equivalently, $F$~is
  visible if $p + \lambda (x-p) \notin P$ for all points $p \in F$ and real
  numbers $\lambda > 0$.  We denote the set of visible faces by~$\Vis(x)$; its
  complement $ \Invis(x) := F(P)_0^1 \setminus \Vis(x) $ is the set of {\it
    invisible\/} faces.
\end{definition}

\begin{figure}[ht]
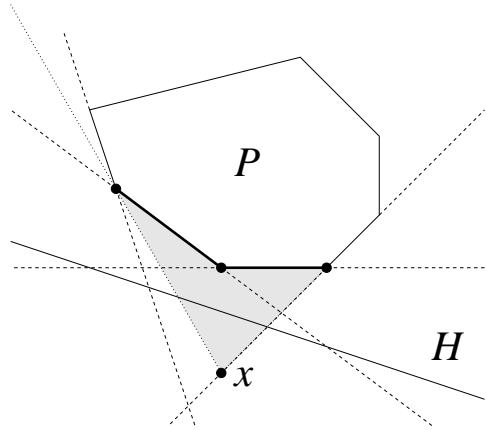

  \begin{center}
    \input visible.pstex_t
    \caption{Visible faces}
    \label{fig:visible}
  \end{center}
\end{figure}

\begin{lemma}
\label{VisibleFacets}
  A facet~$F$ of~$P$ is visible from~$x$ if and only if
  $x$~and $\mathrm {int}\, P$ are on different sides of the affine hyperplane
  spanned by~$F$. A proper non-empty face of~$P$ is visible if and only if
  it is contained in a visible facet of~$P$.
\qed
\end{lemma}

In particular, the sets~$\Vis(x)$ and~$\Invis(x)$ are non-empty.  Since a face of a
visible face is visible itself, $\Vis(x)$ is a subcomplex
while $\Invis(x)$~is an order filter.

\begin{proposition}
\label{VandIareBalls}
  The space $|\Vis(x)|$ is homeomorphic to an $(n-1$)-ball.
  In particular, $\chi (\Vis(x)) = 1$.
\end{proposition}

\begin{proof}%
Applying a translation if necessary we may assume~$x = 0$. Let~$H$
be any hyperplane separating~$0$ and~$P$ (Fig.~\ref{fig:visible}). Let~$C$ denote the cone (with
apex~$0$) on~$P$. Then $C$~is a pointed polyhedral cone, hence $C \cap H$ is a
ball \cite[Theorem~V.1.1]{Ewald-CCAG}. Projection along~$C$ provides a
homeomorphism $|\Vis(x)| \iso C \cap H$.
\end{proof}

\bigbreak

\subsection*{Front and back faces}

\begin{definition}
\label{DefFrontBack}
  A face $F \in F(P)_0^1$ is called a {\it back face\/} with respect to the
  point $x \in \bR^n \setminus \mathrm {int}\,P$
  if for all points $p \in F$ and all real numbers $\lambda > 0$ we have $p +
  \lambda(p-x) \notin P$. The set of back faces is denoted by~$\Back(x)$; its
  complement $\Front(x):= F(P)_0^1 \setminus \Back(x)$ is the set of {\it
    front faces}.
\end{definition}

\begin{figure}[h]
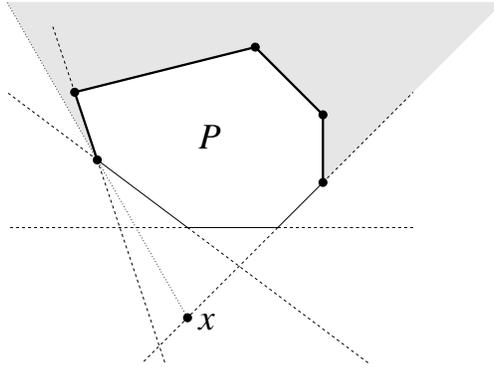

  \begin{center}
    \input back.pstex_t
    \caption{Back faces}
    \label{fig:back}
  \end{center}
\end{figure}

\begin{lemma}
\label{LemmaFrontBack}
  Suppose~$F$ is a facet of~$P$. Then~$F$ is a back face with respect to~$x$
  if and only if~$x$ and~$\mathrm {int}\, P$ are on the same side of the
  affine hyperplane spanned by~$F$. A proper non-empty face~$F$ of~$P$ is a
  back face if and only if it is contained in a facet of~$P$ which is a back face.
\qed
\end{lemma}

In particular, the sets~$\Back(x)$ and~$\Front(x)$ are non-empty.  Since a
face of a back face is a back face itself, $\Back(x)$
is a subcomplex while $\Front(x)$~is an order filter.

By arguments similar to the ones used for the case of visible faces, we can show:

\begin{proposition}
\label{FandBareBalls}
  The space $|\Back(x)|$ is homeomorphic to an $(n-1)$-ball. In
  particular, $\chi (\Back(x)) = 1$. \qed
\end{proposition}

\subsection*{Upper and lower faces}

\begin{definition}
\label{LowUpDef}
  A face $F \in F(P)_0^1$ is called a {\it lower face\/} with respect to the
  direction $x \in \bR^n \setminus \{0\}$
  if for all points $p \in F$ and all real numbers $\lambda > 0$ we have $p -
  \lambda x \notin P$. The set of lower faces is denoted by~$\Low(x)$; its
  complement $\Up(x):= F(P)_0^1 \setminus \Low(x)$ is the set of {\it upper
    faces}.
\end{definition}

\begin{figure}[h]
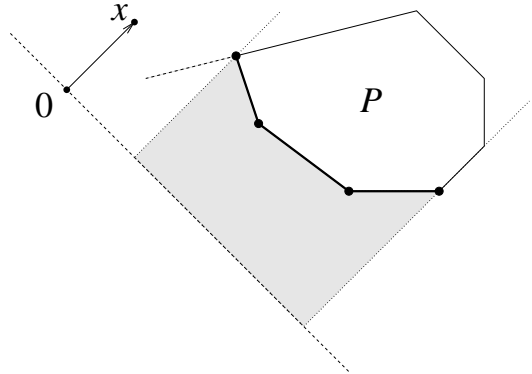

  \begin{center}
    \input lower.pstex_t
    \caption{Lower faces}
    \label{fig:lower}
  \end{center}
\end{figure}

\begin{lemma}
  Suppose~$F$ is a facet of~$P$ with inward pointing normal
  vector~$v$. Then~$F$ is a lower face with respect to~$x$ if and only if
  $\langle x, \, v \rangle > 0$. A proper non-empty face of~$P$ is a lower
  face if and only if it is contained in a facet of~$P$ which is a lower face.
\qed
\end{lemma}

In particular, the sets~$\Low(x)$ and~$\Up(x)$ are non-empty.  Since a face of
a lower face is a lower face itself, $\Low(x)$ is a subcomplex while
$\Up(x)$~is an order filter.

By arguments similar to the ones used for the case of visible faces, we can show:

\begin{proposition}
\label{UandLareBalls}
  The space $|\Low(x)|$ is homeomorphic to an $(n-1)$-ball. In
  particular, $\chi (\Low(x)) = 1$. \qed
\end{proposition}


\section{Barrier cones, tangent cones, and the
  Brian\-chon-Gram theorem}

\label{sec:barr-cones}

Let $P \subset \bR^n$ be a polytope
with non-empty interior. Given a non-empty face $F$
of~$P$ we define the {\it barrier cone\/} $C_F$ of~$P$ at~$F$ as the set of
finite linear combination with non-negative real coefficients spanned by the
set
\[P - F := \{ p - f \,|\, p \in P \hbox{ and } f \in F\} \ .\]
Clearly $C_F$ contains the vector space spanned by~$F-F$ which is the
vector space associated to the affine span of~$F$. This definition
generalises the previous one if $F$ is a vertex of~$P$.

Let $F$ be a non-empty face of~$P$. One should think of the translated
cone $F + C_F = \{f + x \,|\, x \in C_F,\ f \in F\}$ as the cone $C_F$
attached to the face $F$.

For a non-empty proper face $F$ of~$P$ let $T_F$ denote the {\it supporting
  cone\/} (or tangent cone) of~$F$; it is the intersection of all supporting
half-spaces containing $F$ in their boundary. (Of course it is enough to
restrict to facet-defining half-spaces.) By convention $T_P = \bR^n$. Using
\textsc {Farkas\/}' lemma (\cite[\S1.4]{Ziegler-Polytopes} or
\cite[Lemma~I.3.5]{Ewald-CCAG}) it can be shown that $F + C_F = T_F$.
Moreover, every polytope is the intersection of all its supporting
half-spaces, thus $P = \bigcap_{F \in F(P)_0^1} T_F$.

Also of interest are the cones $-F + C_F = \{ -f + x \,|\, x \in C_F,\
f \in F\}$. Up to a reflection at the origin, they can be thought of
as the negatives of the barrier cones, attached to the corresponding
face (\ie, with cones pointing towards the outside of~$P$).

The following theorem is the heart of this paper; expressed in
combinatorial terms, it uses the \textsc{Euler} characteristic to give
specific inclusion-exclusion formul\ae{} for lattice point in (the
interior of)~$P$. Part~(\ref{eq:P}) is known as the
\textsc{Brianchon}-\textsc{Gram} theorem, the remaining two equations
are variations of the theme.

\begin{theorem}
  \label{thm:brion_powerseries}
  Let $P \subset \bR^n$ be an arbitrary $n$-dimensional polytope.
  There are equalities of formal \textsc{Laurent} series
  \begin{eqnarray}
    \label{eq:P} \sum_{F \in F(P)_0} (-1)^{\dim F}
    \mathcal{S}(F+C_F) & = & \mathcal{S}(P) \ ,\\
    \noalign{\vglue 1 em}
    \label{eq:1} \sum_{F \in F(P)_0} (-1)^{\dim F}
    \mathcal{S}(C_F) & = & 1 \ ,\\
    \noalign{\vglue 1 em}
    \label{eq:intP} \sum_{F \in F(P)_0} (-1)^{\dim F}
    \mathcal{S}(-F+C_F) & = & (-1)^n \cdot
    \mathcal{S}(\mathrm{int}\,-P) \ .
  \end{eqnarray}
\end{theorem}

\medskip

\begin{proof}%
We verify Equation~(\ref{eq:P}) first. Fix a vector $\mathbf{a} \in
\bZ^n$. We have to show that the coefficient of~$x^{\mathbf{a}}$ is
the same on both sides of the equation.

If $\mathbf{a} \in P$ then the monomial $x^{\mathbf{a}}$ occurs with
coefficient~$1$ in all the \textsc{Laurent} series
$\mathcal{S}(F+C_F)$ on the left. Thus the coefficient
of~$x^{\mathbf{a}}$ in the sum is the \textsc{Euler} characteristic
of~$P$, which is known to be~$1$. Hence the coefficients
of~$x^{\mathbf{a}}$ agree on the left and right side in this case.

Now assume $\mathbf{a} \notin P$.
Let~$F$ denote a proper non-empty face of~$P$. From
Lemma~\ref{VisibleFacets} and the definition of supporting cones we conclude that
$\mathbf{a} \in T_F = F+C_F$ if and only if $F$~is invisible from~$\mathbf{a}$.
In particular, the coefficient of~$x^{\mathbf{a}}$ in
$\mathcal{S}(F+C_F)$ is~$1$ if $F \in \Invis(\mathbf{a})$, and it is~$0$
if $F \notin \Invis(\mathbf{a})$. In total, the coefficient
of~$x^{\mathbf{a}}$ on the left is
\[\ell = (-1)^n + \sum_{F \in \Invis (\mathbf{a})} (-1)^{\dim(F)} = (-1)^n +
\chi (\Invis(\mathbf{a})) \ ,\]
the extra $(-1)^n$ corresponding to the contribution coming from~$P$. Now by
definition of the \textsc{Euler} characteristic, we have
\[1 = \chi (P) = \chi (\Vis (\mathbf{a})) + \ell \ .\]
Since $|\Vis(\mathbf{a})|$ is a ball by Proposition~\ref{VandIareBalls}, we
infer that $\ell = 0$. Consequently, the monomial~$x^{\mathbf{a}}$ does not
occur on either side of Equation~(\ref{eq:P}), as required.

\smallskip

Next we deal with Equation~(\ref{eq:1}). Observe first that $0 \in C_F$ for
all $F \in F(P)_0$, so the coefficient of $1 = x^0$ is $\chi(P) = 1$.

Now fix any non-zero vector $\mathbf{a} \in \bZ^n$. We have to show that the
coefficient of~$x^{\mathbf{a}}$ is trivial. For a given face $F \in F(P)_0$,
let $N_F := C_F^{\vee}$ denote the dual cone of~$C_F$; it is given by
\[N_F = \left\{ v \in \bR^n \,|\, \forall p \in C_F \colon
   \langle p, \, v \rangle \geq 0 \right\} \ .\]
It can be shown that $N_F$~is the cone of inward pointing normal vectors
of~$F$, and that the dual of~$N_F$, given by
\[N_F^\vee := \left\{ p \in \bR^n \,|\, \forall v \in N_F \colon
   \langle v, \, p \rangle \geq 0 \right\} \ ,\]
is the barrier cone~$C_F$ \cite[\S{}I.4 and~\S{}V.2]{Ewald-CCAG}.

Let $U(\mathbf{a})$ denote the poset of all non-empty proper faces~$F$ of~$P$
satisfying $\mathbf{a} \in C_F$. By the above we have equivalences
\[F \in U(\mathbf{a}) \iff \mathbf{a} \in C_F = N_F^\vee
\iff \forall v \in N_F \colon \langle -\mathbf{a}, \, v \rangle \leq 0 \ .\]
This means that $U(\mathbf{a}) = \Up(-\mathbf{a})$ is the set of upper faces
of~$P$ with respect to~$-\mathbf{a}$ in the sense
of~Definition~\ref{LowUpDef}. Hence the coefficient of~$x^\mathbf{a}$ in the
left-hand side of Equation~(\ref{eq:1}) can be rewritten as
\[ (-1)^n + \sum_{F \in U(\mathbf{a})} (-1)^{\dim(F)} = (-1)^n + \chi
(\Up(-\mathbf{a})) = \chi(P) - \chi (\Low (-\mathbf{a})) = 0\]
where we have used additivity of \textsc{Euler} characteristic and
Proposition~\ref{UandLareBalls} as well.

\smallskip

Finally we discuss Equation~(\ref{eq:intP}).  Fix a point~$\mathbf{a} \in
\bZ^n$ and a face $F \in F(P)_0^1$.
Then $\mathbf{a} \notin -F+C_F$ if and only if there is a facet $G \supseteq F$
of~$P$ such that $\mathbf{a}$ and $\mathrm{int}\,(-P)$ are on the same side of the
affine hyperplane spanned by~$-G$. Such a facet certainly exists if $\mathbf{a} \in
\mathrm{int}\,(-P)$. Hence the only summand on the left contributing
to~$x^{\mathbf{a}}$ is the one corresponding to~$P$, giving a coefficient
$(-1)^n$ as required.

If~$\mathbf{a}$ is not in the interior of~$-P$, Lemma~\ref{LemmaFrontBack},
applied to the polytope~$-P$, shows that $\mathbf{a} \in -F+C_F$ if and only if
$-F$~is a front face of~$-P$ in the sense of Definition~\ref{DefFrontBack}. It
follows from Proposition~\ref{FandBareBalls} and additivity of \textsc{Euler}
characteristics that the contribution to~$x^\mathbf{a}$ is
\[(-1)^n + \sum_{-F \in \Front(-P)} (-1)^{\dim(F)} = \chi (P) - \chi (\Back
(-P)) = 0 \ . \qed\]
\end{proof}

\section{Brion's formula}

\label{sec:brions-formula}

From the \textsc{Brianchon}-\textsc{Gram} theorem we deduce
\textsc{Brion} type formul\ae{} by passing to rational functions. We
follow the treatment as exemplified in \cite{Beck:Brion_formulae}.

Write $\Pi$ for the $\bC [x_1^{\pm 1}, x_2^{\pm 1}, \ldots, x_n^{\pm
  1}]$-submodule of $\bC [[x_1^{\pm 1}, x_2^{\pm 1}, \ldots, x_n^{\pm 1}]]$
generated by formal power series $\mathcal{S}(w + C)$ where $w \in \bR^n$ is
an arbitrary vector, and $C$ is a polyhedral rational cone
\begin{equation}
  \label{eq:cone}
  C = \{ \lambda_1 v_1 + \lambda_2 v_2 + \ldots + \lambda_r v_r \,|\,
  \lambda_1, \ldots, \lambda_r \in \bR_{\geq 0} \}
\end{equation}
where $v_1, \ldots, v_r$ are vectors in $\bZ^n$.  For $r=0$ the cone
$C$ degenerates to
the single point~$0$, so~$\Pi$ contains the series $\mathcal{S}(P)$ for any
polytope~$P$. In general, if $C$ does not contain an affine subspace of
positive dimension (\ie, if $C$ is {\it pointed\/}), the series~$\mathcal{S}(w+C)$
represents a rational function, denoted $\sigma (w+C) \in \bC (x_1, x_2, \ldots,
x_n)$.

The following Lemma has been attributed to \textsc{Brion}, see
\cite[Theorem~2.4]{Beck:Brion_formulae} and
\cite[Theorem~2.1~(i)]{Brion-survey}.

\begin{lemma}
  \label{lem:Brion}
  There is a unique $\bC [x_1^{\pm 1}, x_2^{\pm 1}, \ldots, x_n^{\pm
  1}]$-module homomorphism
  \[\phi \colon \Pi \to \bC (x_1, x_2, \ldots, x_n)\]
  such that $\phi (\mathcal{S}(w+C)) = \sigma(w+C)$ for each pointed cone of the
  form~{\rm (\ref{eq:cone})} where $v_1, \ldots, v_r \in \bZ^n$ and $w
  \in \bR^n$.

  Moreover, if $C$ is of the form~{\rm (\ref{eq:cone})}, and $C$
  contains an affine subspace of positive dimension, then $\phi
  (\mathcal{S}(w+C)) = 0$. \qed
\end{lemma}

\medskip

We now come to the proof of Theorem~\ref{thm:brion}.  We treat
Equation~(\ref{eq:brion_1}) only, the other cases being similar. Note
that for all $F \in F(P)_0$ the barrier cone $C_F$ is a rational
polyhedral cone since the facets of~$P$ admit rational normal vectors.
We can thus apply the homomorphism~$\phi$ from Lemma~\ref{lem:Brion}
to Equation~(\ref{eq:1}), Theorem~\ref{thm:brion_powerseries}. The
results follows immediately if one recalls that $C_F$ contains an
affine subspace of positive dimension if and only if $\dim F \geq 1$,
so all summands coming from faces of positive dimension disappear upon
application of~$\phi$.

\subsection*{Concluding remarks}

Visibility subcomplexes can be used to compute higher sheaf cohomology
of certain line bundles on projective toric varieties; the reader will
easily recognise the similarity between the present paper and the
exposition in \cite{H-nonlin_toric}, Appendix
of~\S2.5. \textsc{Brion}'s theorem can be generalised substantially to
include the case of arbitrary torus-invariant line bundles on
complete toric varieties or, formulated in more
combinatorial terms, arbitrary support functions on complete fans
\cite{H-cohomology}.

\raggedright
\bibliographystyle{alpha}

\begin{thebibliography}{Ewa96}

\bibitem[BHS]{Beck:Brion_formulae}
Matthias Beck, Christian Haase, and Frank Sottile.
\newblock {Theorems of Brion, Lawrence, and Varchenko on rational generating
  functions for cones}.
\newblock arXiv:math.CO/0506466.

\bibitem[Bri88]{Brion-latticepoints}
Michel Brion.
\newblock Points entiers dans les poly\`edres convexes.
\newblock {\em Ann. Sci. \'Ecole Norm. Sup. (4)}, 21(4):653--663, 1988.

\bibitem[Bri96]{Brion-survey}
Michel Brion.
\newblock Polytopes convexes entiers.
\newblock {\em Gaz. Math.}, (67):21--42, 1996.

\bibitem[BV97]{Brion-Vergne}
Michel Brion and Mich{\`e}le Vergne.
\newblock Lattice points in simple polytopes.
\newblock {\em J. Amer. Math. Soc.}, 10(2):371--392, 1997.

\bibitem[Ewa96]{Ewald-CCAG}
G{\"u}nter Ewald.
\newblock {\em Combinatorial convexity and algebraic geometry}.
\newblock Springer-Verlag, New York, 1996.

\bibitem[H{\"u}ta]{H-cohomology}
Thomas H{\"u}ttemann.
\newblock A cohomological interpretation of {B}rion's formula.
\newblock arXiv:math.CO/0607464.

\bibitem[H{\"u}tb]{H-nonlin_toric}
Thomas H{\"u}ttemann.
\newblock {{$K$}-Theory of non-linear projective toric varieties}.
\newblock arXiv:math.KT/0508431.

\bibitem[Zie95]{Ziegler-Polytopes}
G{\"u}nter~M. Ziegler.
\newblock {\em Lectures on polytopes}.
\newblock Springer-Verlag, New York, 1995.

\end{thebibliography}

\end{document}